# COUPLING ALL THE LÉVY STOCHASTIC AREAS OF MULTIDIMENSIONAL BROWNIAN MOTION

By Wilfrid S. Kendall

*University of Warwick*

It is shown how to construct a successful co-adapted coupling of two copies of an $n$-dimensional Brownian motion $(B_1, \ldots, B_n)$ while simultaneously coupling all corresponding copies of Lévy stochastic areas $\int B_i \, dB_j - \int B_j \, dB_i$. It is conjectured that successful co-adapted couplings still exist when the Lévy stochastic areas are replaced by a finite set of multiply iterated path- and time-integrals, subject to algebraic compatibility of the initial conditions.

**1. Introduction.** A probabilistic *coupling* of two random processes is a construction of both processes on the same probability space, building in useful dependence between the two processes. This paper discusses couplings of two Markov processes with the same law of evolution, begun at different points, and constructed so as to join together (to *couple*) at some random time; the coupling is said to be *successful* if the two processes couple within finite time almost surely. There are other kinds of couplings relating to monotonicity, or to approximation; successful couplings are useful for probabilistic gradient estimates, for studying the rate of convergence to statistical equilibrium, for relating behavior of random processes to the geometry of the state-space, and (in more developed formulations) as a key component in perfect simulation algorithms. The present paper focuses on a particular question to do with coupling constructions for Euclidean Brownian motions: namely, whether one can couple successfully not only the Brownian motions themselves, but also sets of path functionals. We shall show that one can couple successfully not only two copies of a Brownian motion $(B_1, \ldots, B_n)$, but at the same time all the corresponding pairs of Lévy stochastic areas $\int B_i \, dB_j - \int B_j \, dB_i$ of the two copies. This appears quite remarkable to the author: one is able to couple so much despite controlling only the correlations between the two copies of the Brownian motion.

Extensive treatments of probabilistic coupling can be found in [16] and [22], so a short summary of the relevant history will suffice. Lindvall [15] was the first to consider coupling for Brownian motion; he described the classic *reflection coupling* (couple two $n$-dimensional Brownian motions by making one of them to be the reflection of the other until they meet). This was followed up by Lindvall and Rogers [17], who discussed generalizations to the case of diffusions. There is a significant distinction to be drawn here. It is typically very much easier to find explicit descriptions of couplings when the two processes in question are required both to be *co-adapted* to the same filtration of $\sigma$-algebras, in particular when the driving Brownian motions have increments which are independent of their common past. In [17] (and throughout the present paper) the search is for co-adapted couplings, and therefore stochastic calculus can be used to provide very explicit descriptions.

Ben Arous, Cranston and Kendall [1] were the first to consider the possibility of what one might call *exotic couplings*, in which one seeks to couple co-adaptively and simultaneously certain path functionals as well as the processes. They described co-adapted couplings for the single stochastic area of planar Brownian motion, and also for the time-integral of scalar Brownian motion $\int B\,dt$. Price [20] showed in her thesis how to extend the second case to couple the twice-iterated time-integral $\int\int B\,ds\,dt$, and Kendall and Price [13] used a different method to show the existence of a successful co-adapted coupling for $B$ and any finite set of iterated time-integrals $\int \cdots \int B\,ds\cdots dt$. The present paper continues this theme by extending the first result of [1] to $n$-dimensional Brownian motion and all possible Lévy stochastic areas $\int B_i\,dB_j - B_j\,dB_i$. It now seems reasonable to formulate a general conjecture that successful co-adapted exotic coupling is possible for any finite combination of multiply iterated path- and time-integrals of Brownian motion (for compatible initial conditions), though it is clear that new ideas will be required to make further progress. The theory of Lie group symmetries of stochastic differential equations supports the expectation that resolution of the general conjecture would lead quickly to coupling constructions for wide classes of hypoelliptic diffusions.

At present the main motivation for studying exotic coupling lies in the importance of coupling as a general concept, and the consequent desirability of understanding how far one can go in coupling large sets of path functionals. However, it does seem not unreasonable to hope for future useful interactions with rough path theory [18], where stochastic areas play a central role, and conceivably also for its use in lifting restrictions on the new methods of exact simulation of stochastic differential equations [2].

It should be noted that there is significant theory concerning nonco-adapted couplings. If the co-adapted constraint is lifted, then one may construct *maximal couplings* [7, 8, 19] which couple at the maximum rate permitted by the total variation bound of the coupling inequality. These couplings have strong relationships with potential theory, and will in general be hard to construct (but see the striking results of Rogers [21] on random walks). Hairer [9] used a restricted kind of nonco-adapted coupling at time $\infty$ to study hypoelliptic diffusions, Hayes and Vigoda [10] used finite look-ahead couplings to gain definite improvements on coupling rate in an application to randomized algorithms, while Burdzy and Kendall [4] studied the cost of the co-adapted property. In our case it is a simple matter to demonstrate the possibility of successful nonco-adapted coupling of Lévy stochastic areas as a consequence of Hörmander regularity of the corresponding $n + n(n-1)/2$-dimensional hypoelliptic diffusion. The point of the present paper is to deliver an explicit successful *co-adapted* coupling construction; the existence of this is *not* implied by regularity theory.

The paper is organized as follows. Section 2 addresses some general considerations related to stochastic control, which help to focus the problem on specific coupling strategies and to introduce notation. Section 3 gives a new approach to the two-dimensional problem treated by [1]; this prepares the way for the main results of the paper which are stated and proved in Section 4: namely, that successful co-adapted coupling is achievable for $n$-dimensional Brownian motion and its $n(n-1)/2$ associated Lévy stochastic areas. The concluding Section 5 considers a couple of complementary issues, and formulates a general question concerning coupling of sets of iterated path integrals, for which the answer is conjectured to be in the affirmative.

**2. Coupling, control and convexity.** It is helpful to bear in mind a stochastic control-theoretic perspective for coupling problems concerning co-adapted Brownian motions (see Borkar [3] for a useful survey on stochastic control; Chen [5] elicits some connections between control and coupling, while Jansons and Metcalfe [12] carry out some numerical investigations). As remarked above, a *co-adapted* coupling of two $n$-dimensional Brownian motions $\underline{\mathbf{A}}$ and $\underline{\mathbf{B}}$ means that $\underline{\mathbf{A}}$ and $\underline{\mathbf{B}}$ are both adapted *as Brownian motions* to the same filtration of $\sigma$-algebras $\{\mathfrak{F}_t : t \geq 0\}$ (thus in particular both increments $\underline{\mathbf{A}}_{t+s} - \underline{\mathbf{A}}_s$ and $\underline{\mathbf{B}}_{t+s} - \underline{\mathbf{B}}_s$ are independent of $\mathfrak{F}_t$). The most general co-adapted coupling can be specified using Itô stochastic calculus:

$$(1) \qquad d\underline{\mathbf{A}} = \underline{\underline{\mathbf{J}}}^T \, d\underline{\mathbf{B}} + \underline{\underline{\widetilde{\mathbf{J}}}}^T \, d\underline{\mathbf{C}},$$

where $\underline{\underline{\mathbf{J}}}$, $\underline{\underline{\widetilde{\mathbf{J}}}}$ are predictable $(n \times n)$-matrix-valued processes, and $\underline{\mathbf{C}}$ is a further $n$-dimensional Brownian motion adapted as a Brownian motion to the filtration $\{\mathfrak{F}_t : t \geq 0\}$ and independent of $\underline{\mathbf{B}}$.

Thus the coupling is specified by giving a control in the form of a pair of predictable matrix-valued processes $\underline{\mathbf{J}}, \underline{\widetilde{\mathbf{J}}}$. These must satisfy certain conditions if (1) is indeed to define an $n$-dimensional Brownian motion $\underline{\mathbf{A}}$: multiplying stochastic differentials to obtain differentials of quadratic variation (following [11]), and bearing in mind the independence of $\underline{\mathbf{B}}$ and $\underline{\mathbf{C}}$, it follows that $\underline{\mathbf{A}}$ is a Brownian motion if and only if the following matrix-valued random measure equation is satisfied:

$$\begin{aligned}
\underline{\mathbf{I}}\, dt &= (d\underline{\mathbf{A}}) \times (d\underline{\mathbf{A}})^T \\
&= (\underline{\mathbf{J}}^T (d\underline{\mathbf{B}}) \times (d\underline{\mathbf{B}})^T \underline{\mathbf{J}}) + (\underline{\widetilde{\mathbf{J}}}^T (d\underline{\mathbf{C}}) \times (d\underline{\mathbf{C}})^T \underline{\widetilde{\mathbf{J}}}) \\
&= (\underline{\mathbf{J}}^T \underline{\mathbf{J}} + \underline{\widetilde{\mathbf{J}}}^T \underline{\widetilde{\mathbf{J}}})\, dt,
\end{aligned}$$
(2)

where $\underline{\mathbf{I}}$ is the $n \times n$ identity matrix.

The matrix process $\underline{\mathbf{J}}$ expresses the infinitesimal correlation of $\underline{\mathbf{C}}$ with $\underline{\mathbf{B}}$: from (2) it follows that such matrix processes are characterized by lying (almost) always in the convex compact set defined by

$$\underline{\mathbf{0}} \leq \underline{\mathbf{J}}^T \underline{\mathbf{J}} \leq \underline{\mathbf{I}},$$
(3)

where $\underline{\mathbf{0}}$ is the $(n \times n)$ zero matrix, and the inequalities are to be interpreted using the usual spectrally based partial ordering for symmetric matrices. An application of the Cauchy–Schwarz inequality to $\underline{\mathbf{v}}^T \underline{\mathbf{J}}^T \underline{\mathbf{J}} \underline{\mathbf{v}}$ shows that the set of extreme points of this compact convex set can be identified as the topological group of orthogonal matrices $O(n)$.

Our coupling problem will be solved by designing a predictable process $\underline{\mathbf{J}}$ such that $\underline{\mathbf{A}}$ and $\underline{\mathbf{B}}$ couple at some finite random time simultaneously with all their stochastic areas ($\int A_i\, dA_j - \int A_j\, dA_i$ coupling with $\int B_i\, dB_j - \int B_j\, dB_i$, etc.). Suppose that it is possible to arrange this in terms of a stochastic control problem which is regular enough to possess an objective function leading to a bounded value function $V(t, \underline{\mathbf{A}}, \underline{\mathbf{B}})$ (where perhaps $t$ is replaced by some other additive functional such as time spent in a specified region). Being a value function, $V(t, \underline{\mathbf{A}}, \underline{\mathbf{B}})$ is a supermartingale in general and is a martingale exactly, when the control $\underline{\mathbf{J}}$ is optimal. If $V(t, \underline{\mathbf{A}}, \underline{\mathbf{B}})$ is appropriately smooth, then Itô's formula may be applied. This together with the Brownian nature of $\underline{\mathbf{A}}$ and $\underline{\mathbf{B}}$ shows

$$\begin{aligned}
dV(t, \underline{\mathbf{A}}, \underline{\mathbf{B}}) = {} & V_0\, dt + \underline{\mathbf{V}}_1^T\, d\underline{\mathbf{A}} + \underline{\mathbf{V}}_2^T\, d\underline{\mathbf{B}} \\
& + \tfrac{1}{2} \operatorname{tr}(\underline{\mathbf{V}}_{11})\, dt + \tfrac{1}{2} \operatorname{tr}(\underline{\mathbf{V}}_{22})\, dt + \operatorname{tr}(\underline{\mathbf{J}} \underline{\mathbf{V}}_{12})\, dt
\end{aligned}$$

(for a fixed orthonormal basis $v_1, \ldots, v_n$, and first- and second-order derivatives $V_0, \underline{\mathbf{V}}_1, \underline{\mathbf{V}}_2, \underline{\mathbf{V}}_{11}, \underline{\mathbf{V}}_{12}, \underline{\mathbf{V}}_{22}$ with dependence on $t, \underline{\mathbf{A}}$ and $\underline{\mathbf{B}}$ suppressed). Thus processes $\underline{\mathbf{J}}$ which are optimal controls for such a regular problem must maximize

$$\operatorname{tr}(\underline{\mathbf{J}} \underline{\mathbf{V}}_{12}),$$

which is linear in $\underline{\underline{\mathbf{J}}}$. It follows that smoothness of an appropriate value function implies that optimal control processes $\underline{\underline{\mathbf{J}}}$ must (almost surely, for almost all time) lie in the region of extreme points of the convex compact region of controls, and so must satisfy the orthogonality condition

$$\underline{\underline{\mathbf{J}}}^T \underline{\underline{\mathbf{J}}} = \underline{\underline{\mathbf{I}}}; \tag{4}$$

in brief, almost surely $\underline{\underline{\mathbf{J}}}(t) \in O(n)$ for almost all times $t$ (and hence $\underline{\underline{\widetilde{\mathbf{J}}}} = \underline{\underline{\mathbf{0}}}$).

The impact of these considerations for our coupling problem is entirely heuristic, since we do not have any particular objective function in mind other than desiring to show that it is possible to couple Brownian motions together with their stochastic areas. (Indeed we will not even check that our resulting coupling strategy is *admissible*, in the sense of being optimal for some objective function: it is not a priori at all clear whether a successful coupling exists and therefore optimality with respect to some arbitrary objective function is of less value than conceptual simplicity!) The above remarks encourage a search for simple couplings among those which use $O(n)$-valued processes $\underline{\underline{\mathbf{J}}}$ to construct $\underline{\mathbf{A}} = \int \underline{\underline{\mathbf{J}}}^T \, d\underline{\mathbf{B}}$ in terms of $\underline{\mathbf{B}}$, without any need of further randomness from $\underline{\mathbf{C}}$. Since $O(n)$ has two topological components, made up of $SO(n)$ and the coset of rotated reflections, it also follows that we should expect to consider coupling strategies which involve discontinuous transitions between one control and another; and this is indeed what may be observed for the successful coupling strategies described in Theorem 4 (for the planar case) and Corollary 7 (for the general case) below.

**3. The planar case.** We first review the planar case (dimension $n = 2$), which permits a simpler treatment than the general $n$-dimensional case but introduces most of the key ideas. The planar case was first dealt with in [1], using controls $\underline{\underline{\mathbf{J}}}$ not all lying in $O(2)$, though it was noted in passing that there was a possibility of coupling using only reflection and synchronous coupling (as defined in Definitions 2 and 3, $\underline{\underline{\mathbf{J}}}$ is a reflection matrix or is an identity matrix). Ben Arous and Lyons have shown in unpublished work how to implement reflection/synchronous coupling for the planar case in a rather direct way, which resembles the low-dimensional case of the Ben Arous et al./Kendall and Price treatment of Kolmogorov diffusions (time-integrals and twice-iterated time-integrals together with scalar Brownian motion). Here we show how reflection/synchronous coupling may be set up using simple and largely state-dependent coupling rules.

First recall from [1] that it is sufficient to couple $(B_1, B_2)$ and $(A_1, A_2)$ together with the *invariant* difference of their stochastic areas,

$$\mathfrak{A} = \int (A_1 \, dA_2 - A_2 \, dA_1) - \int (B_1 \, dB_2 - B_2 \, dB_1) + A_1 B_2 - A_2 B_1. \tag{5}$$

In fact $\mathfrak{A}$ then has a geometric interpretation: it measures the stochastic area swept out by moving first along the $A$ path, then linearly from the end of the $A$ path to the end of the $B$ path, and then back along the $B$ path to its starting point. It turns out to be natural to think of $\mathfrak{A}$ as the $(1,2)$ coordinate of an antisymmetric matrix

$$\underline{\underline{\mathfrak{A}}} = \begin{pmatrix} 0 & \mathfrak{A} \\ -\mathfrak{A} & 0 \end{pmatrix}.$$

Consider the summary quantities

(6)
$$V = \sqrt{(A_1 - B_1)^2 + (A_2 - B_2)^2},$$
$$U = \operatorname{sgn}(\mathfrak{A})\sqrt{\operatorname{tr}(\underline{\underline{\mathfrak{A}}}^T \underline{\underline{\mathfrak{A}}})} = \sqrt{2}\mathfrak{A}.$$

These are semimartingales at least until one of them vanishes. Stochastic calculus can therefore be used to compute the stochastic differential drifts $\operatorname{Drift} dU$ and $\operatorname{Drift} dV$ (the differentials of the locally bounded variation components of the Doob–Meyer semimartingale decompositions of $U$, $V$) and the products of differentials $(dU)^2$, $(dV)^2$ and $(dV) \times (dU)$ (the differentials of the corresponding quadratic variation and covariation processes). In doing this, it is convenient to define the quantities $S_{11}$, $S_{22}$ and $A_{12}$ from the symmetrization and the antisymmetrization of the control $\underline{\mathbf{J}}$: working in orthonormal coordinates based on the vector $\underline{\mathbf{A}} - \underline{\mathbf{B}}$ and its perpendicular,

(7) $\quad \frac{1}{2}(\underline{\underline{\mathbf{J}}} + \underline{\underline{\mathbf{J}}}^T) = \begin{pmatrix} S_{11} & S_{12} \\ S_{21} & S_{22} \end{pmatrix}, \qquad \frac{1}{2}(\underline{\underline{\mathbf{J}}} - \underline{\underline{\mathbf{J}}}^T) = \begin{pmatrix} 0 & A_{12} \\ -A_{12} & 0 \end{pmatrix}.$

The results of these computations are summarized in the following lemma:

LEMMA 1.

$$(dV)^2 = 2(1 - S_{11})\, dt, \qquad \operatorname{Drift} dV = \frac{(1 - S_{22})}{V}\, dt,$$

(8) $\quad (dV) \times (dU) = -2\sqrt{2} A_{12} V\, dt,$

$$(dU)^2 = 4(1 + S_{22})V^2\, dt, \qquad \operatorname{Drift} dU = 2\sqrt{2} A_{12}\, dt.$$

Details of the calculations are left as an exercise for the reader, who may alternatively view them as a special case of the multidimensional case treated in detail in Lemma 5.

Here are two important coupling strategies, defined by specifying the corresponding control $\underline{\underline{\mathbf{J}}}$.

DEFINITION 2. Reflection coupling is defined by choosing $\underline{\underline{\mathbf{J}}}$ to be the orthogonal matrix giving reflection in the line normal to the vector $\underline{\mathbf{A}} - \underline{\mathbf{B}}$; thus $S_{11} = -1$, $S_{22} = 1$, $A_{12} = 0$ in our chosen coordinate system.

Using Lemma 1, reflection coupling yields

$$(dV)^2 = 4\,dt, \qquad \text{Drift}\,dV = 0,$$
$$(9) \qquad (dV) \times (dU) = 0,$$
$$(dU)^2 = 8V^2\,dt, \qquad \text{Drift}\,dU = 0,$$

so that $V$ moves as a scalar Brownian motion at least until it hits 0, and $U$ moves as a scalar Brownian motion subject to a $V$-dependent time-change.

DEFINITION 3. *Synchronous coupling is defined by choosing $\underline{\mathbf{J}}$ to be the identity matrix; thus $S_{11} = S_{22} = 1$, $A_{12} = 0$.*

Using Lemma 1, synchronous coupling yields

$$(dV)^2 = 0, \qquad \text{Drift}\,dV = 0,$$
$$(10) \qquad (dV) \times (dU) = 0,$$
$$(dU)^2 = 8V^2\,dt, \qquad \text{Drift}\,dU = 0,$$

so that $V$ is held constant, while $U$ continues to move as a scalar Brownian motion with rate dependent on $V$ in the same way as for reflection coupling.

Under both these strategies $U$ and $V$ remain semimartingales for all time.

It is possible to derive these results for both couplings without making explicit use of stochastic calculus, simply by considering the geometry of the planar Brownian paths and their invariant difference of areas.

The considerations of Section 2 suggest that if coupling is at all possible for the planar case using only symmetric $\underline{\mathbf{J}}$, then it should be achievable by combining these two controls, since (8) shows that the other two extreme controls ($S_{11} = \pm 1$, $S_{22} = -1$) lead to positive drifts for $V$ without apparent gains for $U$.

Since $U$ scales as $V^2$, and since it is desirable for coupling purposes to reduce the size of $U$ if ever it gets large relative to $V$, it is natural to consider coupling strategies described loosely as follows: for fixed $\kappa > 0$, while $U^2 < \kappa^2 V^4$,    use reflection coupling;
while $U^2 \geq \kappa^2 V^4$,    use synchronous coupling.

This involves a discontinuous change of regime as $(U, V)$ crosses over the boundary $U^2 = \kappa^2 V^4$. The discussion in Section 2 has prepared us to expect such discontinuities. A precise description of a successful strategy of this kind is formulated in the following theorem, which is the principal result of this section.

THEOREM 4. *Suppose that initially $U_0 = 0$ but $V_0 > 0$ (this can always be arranged by first using reflection coupling to make $V$ positive, and then* *using a session of synchronous coupling to reduce $U$ to zero). Fix a small $\varepsilon > 0$; consider the control which alternates between reflection and synchronous couplings using "down-crossings":*

- *if $U^2/V^4$ has not yet visited $\kappa^2$, then use reflection coupling;*
- *if $U^2/V^4$ has attained the level $(\kappa - \varepsilon)^2$ since most recently visiting $\kappa^2$, then use reflection coupling;*
- *otherwise use synchronous coupling.*

*This coupling is almost surely successful in finite time: $(U,V)$ visits $(0,0)$ in finite time.*

Clearly one could consider the limiting case $\varepsilon \to 0$ and use local time and excursion theory; however, it turns out to be simpler to analyze the process as given.

PROOF OF THEOREM 4. Define the indicator process $N^{(\varepsilon)}$ by $N^{(\varepsilon)} = 1$ when either $U^2/V^4$ has not yet visited $\kappa^2$, or $U^2/V^4$ has attained the level $(\kappa - \varepsilon)^2$ since most recently visiting $\kappa^2$, and otherwise set $N^{(\varepsilon)} = 0$. Then the coupling strategy prescribed in the theorem statement corresponds to the stochastic differential system

$$(11) \qquad \begin{aligned} (dV)^2 &= 4N^{(\varepsilon)}\,dt, \qquad \text{Drift}\,dV = 0, \\ (dV) &\times (dU) = 0, \\ (dU)^2 &= 8V^2\,dt, \qquad \text{Drift}\,dU = 0. \end{aligned}$$

This is solvable up to the time when $U$ and $V$ both vanish: one may piece together solutions of the smooth systems defined by (9) and (10). Under this stochastic differential system the process $V$ evolves as a Brownian motion of rate 4 interrupted only when $U^2/V^4$ makes down-crossings from $\kappa^2$ to $(\kappa - \varepsilon)^2$, and during these interruptions $V$ is frozen. These down-crossings each take a finite amount of time, and only finitely many occur in bounded closed time intervals before $U$ and $V$ both vanish; consequently $V$ either hits 0 at a finite time or converges to 0 at time $\infty$. Since $V$ is constant when $U^2/V^4 \geq \kappa^2$, continuity considerations show that $U/V^2 \to 0$ as $V \to 0$, and therefore coupling must occur when $V$ hits 0. Thus the crux of the matter is, will $V \to 0$ at a finite time?

To analyze this question, apply Lamperti's [14] observation (as used to great effect in [23], e.g.) to the stochastic differential system (11). Consider a random time-change under which $K = \log(V)$ behaves as an (interrupted) Brownian motion with constant negative drift. The time-change $\tau(t)$ is defined by

$$(12) \qquad 4\,dt = V^2\,d\tau.$$

Writing $W = U/V^2$, the stochastic system for $K$ and $W$ then follows by Itô's formula:

$$(dK)^2 = N^{(\varepsilon)} d\tau, \qquad \text{Drift } dK = -\tfrac{1}{2} N^{(\varepsilon)} d\tau,$$

$$(13) \quad (dK) \times (dW) = 2 N^{(\varepsilon)} W \, d\tau,$$

$$(dW)^2 = 2(1 + 2 N^{(\varepsilon)} W^2) \, d\tau, \qquad \text{Drift } dW = 3 N^{(\varepsilon)} W \, d\tau.$$

It is required to show that elapsed $t$-time until $K \to -\infty$ (equivalently $V = 0$) is finite, which is equivalent to showing

$$(14) \qquad \int_0^\infty e^{2K} \, d\tau < \infty.$$

Since $V$ diffuses as Brownian motion of rate 4 when $N^{(\varepsilon)} = 1$ and is otherwise frozen, it follows that the integral $\int_0^\infty N^{(\varepsilon)} e^{2K} \, d\tau$ is a Brownian first-passage time and therefore is finite. Accordingly, it is enough to show

$$(15) \qquad \int_0^\infty (1 - N^{(\varepsilon)}) e^{2K} \, d\tau < \infty.$$

Let $\sigma_n^s < \sigma_n^f$ be the start and finish times (in $\tau$-time-scale) of the $n$th down-crossing of $W^2 = U^2/V^4$ from $\kappa^2$ to $(\kappa - \varepsilon)^2$. But $N_\tau^{(\varepsilon)} = 0$ exactly when $\tau$ lies in the union of the stopping-time intervals $[\sigma_n^s, \sigma_n^f]$, so therefore

$$(16) \qquad \int_0^\infty (1 - N^{(\varepsilon)}) e^{2K} \, d\tau = \sum_{n=1}^\infty e^{2 K_{\sigma_n^s}} (\sigma_n^f - \sigma_n^s),$$

since $V = e^K$ remains constant for $\tau \in [\sigma_n^s, \sigma_n^f]$.

Conditional on $K_{\sigma_n^s} : n = 1, 2, \ldots$, the durations $\sigma_n^f - \sigma_n^s$ are independent Brownian first-passage times of different rates. Consequently

$$(17) \quad \begin{aligned} & \mathbb{E}\left[ \exp\left( -\sum_{n=1}^\infty e^{2 K_{\sigma_n^s}} (\sigma_n^f - \sigma_n^s) \right) \Big| K_{\sigma_n^s} : n = 1, 2, \ldots \right] \\ & = \exp\left( -\sum_{n=1}^\infty e^{K_{\sigma_n^s}} \times \varepsilon \right), \end{aligned}$$

using the formula for the moment-generating function of a Brownian first-passage time.

Consider now the times $\sigma_2^s - \sigma_1^f$, $\sigma_3^s - \sigma_2^f, \ldots$ between successive down-crossings. These are independent, identically distributed, and of finite mean, since their common distribution is the $\tau$-time for the regular real-line diffusion $W$ (with $N^{(\varepsilon)} = 1$) to hit one of $\pm\kappa$ when started at $\kappa - \varepsilon$. Thus by the strong law of large numbers it follows that almost surely

$$\frac{1}{n} \sum_{m=1}^n (\sigma_m^s - \sigma_{m-1}^f) \to \mathbb{E}[\sigma_2^s - \sigma_1^f] > 0$$

(defining $\sigma_0^f = 0$).

But equally $K$ is a Brownian motion with constant drift of $-\frac{1}{2}$ on the interrupted $\tau$-time-scale $\int N^{(\varepsilon)} \, d\tau$, and therefore almost surely

$$\frac{K_{\sigma_n^s}}{\int_0^{\sigma_n^s} N^{(\varepsilon)} \, d\tau} = \frac{K_{\sigma_n^s}}{\sum_{m=1}^n (\sigma_m^s - \sigma_{m-1}^f)} \to -\frac{1}{2}.$$

It follows that almost surely

(18) $$\frac{K_{\sigma_n^s}}{n} \to -\frac{1}{2} \mathbb{E}[\sigma_2^s - \sigma_1^f] < 0,$$

and hence $\sum_{n=1}^\infty e^{K_{\sigma_n^s}}$ is almost surely finite.

Consequently (17) shows that

$$\mathbb{E}\left[\exp\left(-\sum_{n=1}^\infty e^{2K_{\sigma_n^s}}(\sigma_n^f - \sigma_n^s)\right) \bigg| K_{\sigma_n^s} \colon n = 1, 2, \ldots\right]$$

is almost surely positive, and so

$$\sum_{n=1}^\infty e^{2K_{\sigma_n^s}}(\sigma_n^f - \sigma_n^s)$$

has a positive chance of being finite, even when conditioned on $K_{\sigma_n^s} \colon n = 1, 2, \ldots$. But the $e^{2K_{\sigma_n^s}}(\sigma_n^f - \sigma_n^s)$ are independent under this conditioning, and so by the Kolmogorov zero-one law and (16) it follows that

(19) $$\int_0^\infty (1 - N^{(\varepsilon)}) e^{2K} \, d\tau < \infty$$

with probability 1. It follows that coupling under this strategy almost surely succeeds after a finite time. □

Further development of this line of reasoning delivers an explicit construction of the limiting case $\varepsilon \to 0$ using local time and excursion theory, a single elliptic partial differential equation for the moment-generating function

$$\mathbb{E}[\exp(-pT)]$$

of the coupling time $T$ for all $p$ using scaling, and estimates for exceedance probabilities of the coupling time. We do not consider these topics here, but instead proceed to the multidimensional case.

**4. The $n$-dimensional case.** The first step is to establish the stochastic differential system (6) for Euclidean separation and invariant difference of stochastic areas, working in general $n$-space $(n > 2)$. First introduce new coordinates based on $\underline{\mathbf{X}} = \underline{\mathbf{A}} - \underline{\mathbf{B}}$ and $\underline{\mathbf{Y}} = \underline{\mathbf{A}} + \underline{\mathbf{B}}$, where $A$ and $B$ are co-adapted $n$-dimensional Brownian motions satisfying (1). Using Itô calculus for the vectors $d\underline{\mathbf{X}}$ and $d\underline{\mathbf{Y}}$,

$$
\begin{aligned}
d\underline{\mathbf{X}}\, d\underline{\mathbf{X}}^T &= 2(\underline{\mathbf{I}} - \underline{\mathbf{S}})\, dt, \qquad \text{Drift}\, d\underline{\mathbf{X}} = \underline{\mathbf{0}}, \\
d\underline{\mathbf{Y}}\, d\underline{\mathbf{X}}^T &= 2\underline{\mathbf{A}}\, dt, \\
d\underline{\mathbf{Y}}\, d\underline{\mathbf{Y}}^T &= 2(\underline{\mathbf{I}} + \underline{\mathbf{S}})\, dt, \qquad \text{Drift}\, d\underline{\mathbf{Y}} = \underline{\mathbf{0}},
\end{aligned}
\tag{20}
$$

where $\underline{\mathbf{S}} = \tfrac{1}{2}(\underline{\mathbf{J}} + \underline{\mathbf{J}}^T)$ and $\underline{\mathbf{A}} = \tfrac{1}{2}(\underline{\mathbf{J}} - \underline{\mathbf{J}}^T)$ are the symmetrized and antisymmetrized matrices corresponding to $\underline{\mathbf{J}}$.

Applying the Itô formula to $V^2 = \underline{\mathbf{X}}^T \underline{\mathbf{X}}$ (the square of the length of $\underline{\mathbf{X}}$), it follows that while $V$ remains positive

$$
(dV)^2 = 2(1 - \underline{\nu}^T \underline{\mathbf{S}}\, \underline{\nu})\, dt, \qquad \text{Drift}\, dV = \frac{n - 1 - (\operatorname{tr} \underline{\mathbf{S}} - \underline{\nu}^T \underline{\mathbf{S}}\, \underline{\nu})}{V}\, dt,
\tag{21}
$$

where $\underline{\nu} = \underline{\mathbf{X}}/V$ is a normalized *configuration vector* defined by $\underline{\mathbf{X}} = \underline{\mathbf{A}} - \underline{\mathbf{B}}$.

Now consider the antisymmetric matrix $\underline{\mathfrak{A}}$ determined by invariant differences of stochastic areas of the form of (5):

$$
\mathfrak{A}_{ij} = \int (A_i\, dA_j - A_j\, dA_i) - \int (B_i\, dB_j - B_j\, dB_i) + A_i B_j - A_j B_i.
$$

Since $\underline{\mathbf{A}} = \tfrac{1}{2}(\underline{\mathbf{Y}} + \underline{\mathbf{X}})$ and $\underline{\mathbf{B}} = \tfrac{1}{2}(\underline{\mathbf{Y}} - \underline{\mathbf{X}})$, calculation shows

$$
d\mathfrak{A}_{ij} = X_i\, dY_j - X_j\, dY_i - 2A_{ij}\, dt.
\tag{22}
$$

Hence

$$
\begin{aligned}
d\mathfrak{A}_{ij} &\times d\mathfrak{A}_{rs} \\
&= X_i X_r\, dY_j\, dY_s - X_j X_r\, dY_i\, dY_s - X_i X_s\, dY_j\, dY_r + X_j X_s\, dY_i\, dY_r \\
&= 2(X_i X_r (\underline{\mathbf{I}} + \underline{\mathbf{S}})_{js} - X_j X_r (\underline{\mathbf{I}} + \underline{\mathbf{S}})_{is} \\
&\quad - X_i X_s (\underline{\mathbf{I}} + \underline{\mathbf{S}})_{jr} + X_j X_s (\underline{\mathbf{I}} + \underline{\mathbf{S}})_{ir})\, dt.
\end{aligned}
\tag{23}
$$

Setting $U = \operatorname{tr}(\underline{\mathfrak{A}}^T \underline{\mathfrak{A}})$, since

$$
d(U^2) = 2U\, dU + (dU)^2 = \sum_i \sum_j (2\mathfrak{A}_{ij}\, d\mathfrak{A}_{ij} + (d\mathfrak{A}_{ij})^2)
\tag{24}
$$

it follows

$$
\begin{aligned}
4U^2 (dU)^2 &= 4 \sum_i \sum_j \sum_r \sum_s \mathfrak{A}_{ij} \mathfrak{A}_{rs}\, d\mathfrak{A}_{ij}\, d\mathfrak{A}_{rs} \\
&= 32 \sum_i \sum_j \sum_r \sum_s \mathfrak{A}_{ij} \mathfrak{A}_{rs} X_i X_r (\underline{\mathbf{I}} + \underline{\mathbf{S}})_{js}\, dt \\
&= 32 \underline{\mathbf{X}}^T \underline{\mathfrak{A}}^T (\underline{\mathbf{I}} + \underline{\mathbf{S}}) \underline{\mathfrak{A}}\, \underline{\mathbf{X}}\, dt = 32 \underline{\nu}^T \underline{\mathfrak{z}}^T (\underline{\mathbf{I}} + \underline{\mathbf{S}}) \underline{\mathfrak{z}}\, \underline{\nu}\, U^2 V^2\, dt.
\end{aligned}
\tag{25}
$$

Here $\underline{\underline{\mathfrak{z}}} = \underline{\underline{\mathfrak{A}}}/U$ is a normalized *configuration matrix* [with $\text{tr}(\underline{\underline{\mathfrak{z}}}^T \underline{\underline{\mathfrak{z}}}) = 1$, antisymmetric so $\underline{\underline{\mathfrak{z}}}^T = -\underline{\underline{\mathfrak{z}}}$ and $\underline{\nu}^T \underline{\underline{\mathfrak{z}}} \underline{\nu} = 0$]. The second line of (25) follows from the first by applying (23) and then exploiting the symmetry of $\underline{\underline{\mathbf{I}}} + \underline{\underline{\mathbf{S}}}$ and the antisymmetry of $\underline{\underline{\mathfrak{A}}}$.

On the other hand, from (24),

$$2U \,\text{Drift}\, dU$$
$$= \text{Drift}\, d(U^2) - (dU)^2$$
$$= \text{Drift} \sum_i \sum_j (2\mathfrak{A}_{ij}\, d\mathfrak{A}_{ij} + (d\mathfrak{A}_{ij})^2) - (dU)^2$$
$$= 4\,\text{tr}(\underline{\underline{\mathfrak{A}}}^T \underline{\underline{\mathbf{A}}})\, dt - (dU)^2$$
$$(26) \qquad + \sum_i \sum_j (X_i^2 (dY_j)^2 + X_j^2 (dY_i)^2 - 2X_i X_j\, dY_j\, dY_i)\, dt$$
$$= 4\,\text{tr}(\underline{\underline{\mathfrak{A}}}^T \underline{\underline{\mathbf{A}}})\, dt - (dU)^2$$
$$+ 2 \sum_i \sum_j (2X_i^2 (\underline{\underline{\mathbf{I}}} + \underline{\underline{\mathbf{S}}})_{jj} - 2X_i X_j (\underline{\underline{\mathbf{I}}} + \underline{\underline{\mathbf{S}}})_{ij})\, dt$$
$$= 4\,\text{tr}(\underline{\underline{\mathfrak{A}}}^T \underline{\underline{\mathbf{A}}})\, dt - (dU)^2 + 4(\text{tr}(\underline{\underline{\mathbf{I}}} + \underline{\underline{\mathbf{S}}}) - \underline{\nu}^T(\underline{\underline{\mathbf{I}}} + \underline{\underline{\mathbf{S}}})\underline{\nu})V^2\, dt$$
$$= 4\,\text{tr}(\underline{\underline{\mathfrak{A}}}^T \underline{\underline{\mathbf{A}}})\, dt + 4(n-1 + \text{tr}\,\underline{\underline{\mathbf{S}}} - \underline{\nu}^T \underline{\underline{\mathbf{S}}} \underline{\nu} - 2\underline{\nu}^T \underline{\underline{\mathfrak{z}}}^T (\underline{\underline{\mathbf{I}}} + \underline{\underline{\mathbf{S}}}) \underline{\underline{\mathfrak{z}}} \underline{\nu})V^2\, dt,$$

where the last line is derived from (25), evaluating $\text{tr}\,\underline{\underline{\mathbf{I}}} = n$, $\underline{\nu}^T \underline{\underline{\mathbf{I}}} \underline{\nu} = 1$.

From (25) and (26) taken together,

$$(dU)^2 = 8\underline{\nu}^T \underline{\underline{\mathfrak{z}}}^T (\underline{\underline{\mathbf{I}}} + \underline{\underline{\mathbf{S}}}) \underline{\underline{\mathfrak{z}}} \underline{\nu} V^2\, dt,$$

$$(27) \quad \text{Drift}\, dU = 2\,\text{tr}(\underline{\underline{\mathfrak{z}}}^T \underline{\underline{\mathbf{A}}})\, dt$$
$$+ 2(n - 1 + \text{tr}\,\underline{\underline{\mathbf{S}}} - \underline{\nu}^T \underline{\underline{\mathbf{S}}} \underline{\nu} - 2\underline{\nu}^T \underline{\underline{\mathfrak{z}}}^T (\underline{\underline{\mathbf{I}}} + \underline{\underline{\mathbf{S}}}) \underline{\underline{\mathfrak{z}}} \underline{\nu}) \frac{V^2}{U}\, dt.$$

Finally, using the antisymmetry of $\underline{\underline{\mathfrak{A}}}$,

$$d(U^2)\, d(V^2) = 4VU\, dV\, dU$$
$$= 4 \sum_i X_i\, dX_i \sum_r \sum_s \mathfrak{A}_{rs}\, d\mathfrak{A}_{rs}$$
$$= -16 \underline{\nu}^T \underline{\underline{\mathfrak{z}}}^T \underline{\underline{\mathbf{A}}} \underline{\nu} UV^2\, dt$$

and so finally

$$(28) \qquad dU\, dV = -4\underline{\nu}^T \underline{\underline{\mathfrak{z}}}^T \underline{\underline{\mathbf{A}}} \underline{\nu} V\, dt.$$

Following the procedure of the planar case, now consider the behavior of $K = \log(V)$. As in Section 3, define $W = U/V^2$; however, we will consider the behavior of $K$ together with that of $H = \log(U)$ rather than that of $W = \exp(H - 2K)$. The next lemma follows from the calculations in this section so far.

LEMMA 5. *For a general control $\underline{\underline{\mathbf{J}}}$ (with symmetric and antisymmetric components $\underline{\underline{\mathbf{S}}}$ and $\underline{\underline{\mathbf{A}}}$), and defining a new ($\tau$-)time-scale by $4\,dt = V^2\,d\tau$ as in Section 3,*

$$(dK)^2 = \tfrac{1}{2}(1 - \underline{\nu}^T\underline{\underline{\mathbf{S}}}\,\underline{\nu})\,d\tau,$$

$$\text{Drift}\,dK = \tfrac{1}{4}(n - \operatorname{tr}\underline{\underline{\mathbf{S}}} - 2(1 - \underline{\nu}^T\underline{\underline{\mathbf{S}}}\,\underline{\nu}))\,d\tau,$$

$$(dK) \times (dH) = -(\underline{\nu}^T\underline{\underline{\mathbf{3}}}^T\underline{\underline{\mathbf{A}}}\,\underline{\nu})\frac{1}{W}\,d\tau,$$

(29)
$$(dH)^2 = 2\underline{\nu}^T\underline{\underline{\mathbf{3}}}^T(\underline{\underline{\mathbf{I}}} + \underline{\underline{\mathbf{S}}})\underline{\underline{\mathbf{3}}}\,\underline{\nu}\frac{1}{W^2}\,d\tau,$$

$$\text{Drift}\,dH = \frac{1}{2}\operatorname{tr}(\underline{\underline{\mathbf{3}}}^T\underline{\underline{\mathbf{A}}})\frac{1}{W}\,d\tau$$
$$+ \frac{1}{2}(n - 1 + \operatorname{tr}\underline{\underline{\mathbf{S}}} - \underline{\nu}^T\underline{\underline{\mathbf{S}}}\,\underline{\nu} - 4\underline{\nu}^T\underline{\underline{\mathbf{3}}}^T(\underline{\underline{\mathbf{I}}} + \underline{\underline{\mathbf{S}}})\underline{\underline{\mathbf{3}}}\,\underline{\nu})\frac{1}{W^2}\,d\tau.$$

PROOF. Use (21), (27) and (28), and Itô's formula. □

The special cases of reflection and synchronous coupling now follow directly. Reflection coupling is defined by

(30) $$\underline{\underline{\mathbf{J}}}^{\text{reflection}} = \underline{\underline{\mathbf{I}}} - 2\underline{\nu}\underline{\nu}^T,$$

which implies

$$\underline{\underline{\mathbf{S}}} = \underline{\underline{\mathbf{J}}}^{\text{reflection}}, \qquad \underline{\underline{\mathbf{A}}} = \underline{\underline{\mathbf{0}}},$$
$$\operatorname{tr}\underline{\underline{\mathbf{S}}} = n - 2, \qquad \underline{\nu}^T\underline{\underline{\mathbf{S}}}\,\underline{\nu} = -1, \qquad \underline{\underline{\mathbf{S}}}\underline{\underline{\mathbf{3}}}\,\underline{\nu} = \underline{\underline{\mathbf{3}}}\,\underline{\nu},$$

and consequently

$$(dK)^2 = d\tau, \qquad \text{Drift}\,dK = -\tfrac{1}{2}\,d\tau,$$

(31) $(dK) \times (dH) = 0,$

$$(dH)^2 = 4\|\underline{\underline{\mathbf{3}}}\,\underline{\nu}\|^2\frac{d\tau}{W^2}, \qquad \text{Drift}\,dH = (n - 1 - 4\|\underline{\underline{\mathbf{3}}}\,\underline{\nu}\|^2)\frac{d\tau}{W^2}.$$

Synchronous coupling is defined by

(32) $$\underline{\underline{\mathbf{J}}}^{\text{synchronous}} = \underline{\underline{\mathbf{I}}},$$

which implies

$$\underline{\underline{\mathbf{S}}} = \underline{\underline{\mathbf{J}}}^{\text{synchronous}}, \qquad \underline{\underline{\mathbf{A}}} = \underline{\underline{\mathbf{0}}},$$
$$\operatorname{tr}\underline{\underline{\mathbf{S}}} = n, \qquad \underline{\nu}^T\underline{\underline{\mathbf{S}}}\,\underline{\nu} = 1, \qquad \underline{\underline{\mathbf{S}}}\underline{\underline{\mathbf{3}}}\,\underline{\nu} = \underline{\underline{\mathbf{3}}}\,\underline{\nu},$$

and consequently

$$(dK)^2 = 0, \qquad \text{Drift}\, dK = 0,$$

(33) $\quad (dK) \times (dH) = 0,$

$$(dH)^2 = 4\|\underline{\underline{\mathfrak{Z}}}\nu\|^2 \frac{d\tau}{W^2}, \qquad \text{Drift}\, dH = (n - 1 - 4\|\underline{\underline{\mathfrak{Z}}}\nu\|^2)\frac{d\tau}{W^2}.$$

Note that $\|\underline{\underline{\mathfrak{Z}}}\nu\|^2$ is bounded above by $1/2$, since the nonzero eigenvalues of the antisymmetric matrix $\underline{\underline{\mathfrak{Z}}}$ all have multiplicity 2 and the sum of squared eigenvalues is $\operatorname{tr}(\underline{\underline{\mathfrak{Z}}}^T\underline{\underline{\mathfrak{Z}}}) = 1$. So if $n \geq 3$, then $H$ is a nonconstant submartingale under both controls; it follows that there is no hope of coupling higher-dimensional stochastic areas by using only synchronous and reflection coupling. Instead we analyze the more complicated case of general orthogonal-matrix controls.

Consider the case of a *rotation coupling* defined adaptively by a matrix exponential

(34) $$\underline{\underline{\mathbf{J}}}^{\text{rotation}}(\theta\underline{\underline{\mathfrak{J}}}) = \exp(\theta\underline{\underline{\mathfrak{J}}}).$$

Here $\underline{\underline{\mathfrak{J}}}$ is an antisymmetric matrix satisfying $\operatorname{tr}(\underline{\underline{\mathfrak{J}}}^T\underline{\underline{\mathfrak{J}}}) = 1$, so that $\underline{\underline{\mathbf{J}}}^{\text{rotation}}(\theta\underline{\underline{\mathfrak{J}}})$ is indeed a rotation matrix, and moreover a finite Taylor series expansion produces an approximation which can be bounded:

$$\underline{\underline{\mathbf{S}}} = \cosh(\theta\underline{\underline{\mathfrak{J}}}) = \tfrac{1}{2}(\underline{\underline{\mathbf{J}}}^{\text{rotation}}(\theta\underline{\underline{\mathfrak{J}}}) + \underline{\underline{\mathbf{J}}}^{\text{rotation}}(-\theta\underline{\underline{\mathfrak{J}}})) = \underline{\underline{\mathbf{I}}} - \frac{\theta^2}{2}\underline{\underline{\mathfrak{J}}}^T\underline{\underline{\mathfrak{J}}} + \theta^4\underline{\underline{\mathbf{O}}}(1),$$

$$\underline{\underline{\mathbf{A}}} = \sinh(\theta\underline{\underline{\mathfrak{J}}}) = \tfrac{1}{2}(\underline{\underline{\mathbf{J}}}^{\text{rotation}}(\theta\underline{\underline{\mathfrak{J}}}) - \underline{\underline{\mathbf{J}}}^{\text{rotation}}(-\theta\underline{\underline{\mathfrak{J}}})) = \theta\underline{\underline{\mathfrak{J}}} + \theta^3\underline{\underline{\mathbf{O}}}(1).$$

Here the $\underline{\underline{\mathbf{O}}}(1)$ terms in the errors signify matrices which vary from line to line but which can be bounded uniformly in $\theta$ and $\underline{\underline{\mathfrak{J}}}$. Hence

$$\operatorname{tr}\underline{\underline{\mathbf{S}}} = n - \frac{\theta^2}{2} + \theta^4 O(1), \qquad \nu^T\underline{\underline{\mathbf{S}}}\nu = 1 - \frac{\theta^2}{2}\|\underline{\underline{\mathfrak{J}}}\nu\|^2 + \theta^4 O(1),$$

$$\operatorname{tr}(\underline{\underline{\mathfrak{Z}}}^T\underline{\underline{\mathbf{A}}}) = \theta\operatorname{tr}(\underline{\underline{\mathfrak{Z}}}^T\underline{\underline{\mathfrak{J}}}) + \theta^3 O(1), \qquad \nu^T\underline{\underline{\mathfrak{Z}}}^T\underline{\underline{\mathbf{A}}}\nu = \theta\langle\underline{\underline{\mathfrak{Z}}}\nu, \underline{\underline{\mathfrak{J}}}\nu\rangle + \theta^3 O(1),$$

$$\nu^T\underline{\underline{\mathfrak{Z}}}^T(\underline{\underline{\mathbf{I}}} + \underline{\underline{\mathbf{S}}})\underline{\underline{\mathfrak{Z}}}\nu = 2\|\underline{\underline{\mathfrak{Z}}}\nu\|^2 + \theta^2 O(1),$$

where again the $O(1)$ terms in the errors (both here and in the following exposition) vary from line to line but are bounded uniformly in $\theta$, $\underline{\underline{\mathfrak{J}}}$ and the configuration matrix $\underline{\underline{\mathfrak{Z}}}$. For the sake of simplicity we choose $\theta = -\gamma/W$, $\underline{\underline{\mathfrak{J}}} = \underline{\underline{\mathfrak{Z}}}$, and consider the effects of applying the adaptive rotational control $\underline{\underline{\mathbf{J}}} = \underline{\underline{\mathbf{J}}}^{\text{rotation}}(-\gamma\underline{\underline{\mathfrak{Z}}}/W)$:

$$(dK)^2 = \frac{\gamma^2}{4}\|\underline{\underline{\mathfrak{Z}}}\nu\|^2\frac{d\tau}{W^2} + \frac{\gamma^4}{W^4}O(1)\,d\tau,$$

$$\text{Drift } dK = \frac{\gamma^2}{8}(1 - 2\|\underline{\underline{\mathfrak{Z}}}\underline{\nu}\|^2)\frac{d\tau}{W^2} + \frac{\gamma^4}{W^4}O(1)\,d\tau,$$

(35)
$$(dK) \times (dH) = \gamma\|\underline{\underline{\mathfrak{Z}}}\underline{\nu}\|^2\frac{d\tau}{W^2} + \frac{\gamma^3}{W^4}O(1)\,d\tau,$$

$$(dH)^2 = 4\|\underline{\underline{\mathfrak{Z}}}\underline{\nu}\|^2\frac{d\tau}{W^2} + \frac{\gamma^2}{W^4}O(1)\,d\tau,$$

$$\text{Drift } dH = -\left(\frac{\gamma}{2} - (n - 1 - 4\|\underline{\underline{\mathfrak{Z}}}\underline{\nu}\|^2)\right)\frac{d\tau}{W^2} + \frac{\gamma^2 + \gamma^3}{W^4}O(1)\,d\tau.$$

The antisymmetric component of the control contributes a crucial $\frac{-\gamma\,d\tau}{2W^2}$ term to the drift of $H$. This can be used to make $H$ a supermartingale. (Incidentally, the choice $\underline{\underline{\mathfrak{J}}} = \underline{\underline{\mathfrak{Z}}}$ maximizes this particular term.)

This motivates a direct construction of a successful coupling strategy, using a mixture of $\underline{\underline{\mathbf{J}}}^{\text{reflection}}$ and $\underline{\underline{\mathbf{J}}}^{\text{rotation}}(-\gamma\underline{\underline{\mathfrak{Z}}}/W)$ with adaptive choices of parameters. This delivers a positive chance of successful coupling for large initial values $W_0$ of $W$:

THEOREM 6. *Consider the adaptively mixed coupling*

$$\underline{\underline{\mathbf{J}}} = \frac{\delta}{W^2}\underline{\underline{\mathbf{J}}}^{\text{reflection}} + \left(1 - \frac{\delta}{W^2}\right)\underline{\underline{\mathbf{J}}}^{\text{rotation}}\left(-\frac{\gamma}{W}\underline{\underline{\mathfrak{Z}}}\right),$$

(36)
$$\delta = \delta(\underline{\underline{\mathfrak{Z}}},\underline{\nu}) = 2\left(\mu_K + \frac{\gamma^2}{8}(1 - 2\|\underline{\underline{\mathfrak{Z}}}\underline{\nu}\|^2)\right),$$

$$\gamma = \gamma(\underline{\underline{\mathfrak{Z}}},\underline{\nu}) = 2(\mu_H + n - 1 - 4\|\underline{\underline{\mathfrak{Z}}}\underline{\nu}\|^2),$$

*defined so long as*

$$W^2 > \delta_0 = 2\mu_K + (\mu_H + n - 1)^2.$$

*This coupling strategy has a positive probability of being successful within finite time if $W_0^2 > w^{(\varepsilon)}$, where $W_0$ is the initial value of $W$ at time $0$ and $w^{(\varepsilon)}$ is a certain finite threshold defined by (39) below, so long as we choose*

(37)
$$0 < \mu_K < \mu_H < \mu_K.$$

*Moreover, the coupling strategy will succeed almost surely if $W$ stays above the threshold $w^{(\varepsilon)}$ for all time.*

Recall from the discussion after (33) that $\|\underline{\underline{\mathfrak{Z}}}\underline{\nu}\|^2$ is bounded above by $1/2$. So $\delta - 2\mu_K$ as given above is always nonnegative (as is $\gamma - 2\mu_H$ if $n \geq 3$).

PROOF OF THEOREM 6. The effect of the mixed control can be evaluated as a convex combination of the systems of reflection coupling (30) and

rotation coupling (35):

$$(dK)^2 = (2\mu_K + (\mu_H + n - 1 - 4\|\underline{\underline{3}}\underline{\nu}\|^2)^2(1 - \|\underline{\underline{3}}\underline{\nu}\|^2))\frac{d\tau}{W^2}$$
$$+ \frac{O(1)}{W^4}d\tau,$$
$$\text{Drift } dK = -\mu_K \frac{d\tau}{W^2} + \frac{O(1)}{W^4}d\tau,$$
(38)
$$(dK) \times (dH) = 2(\mu_H + n - 1 - 4\|\underline{\underline{3}}\underline{\nu}\|^2)\|\underline{\underline{3}}\underline{\nu}\|^2\frac{d\tau}{W^2} + \frac{O(1)}{W^4}d\tau,$$
$$(dH)^2 = 4\|\underline{\underline{3}}\underline{\nu}\|^2\frac{d\tau}{W^2} + \frac{O(1)}{W^4}d\tau,$$
$$\text{Drift } dH = -\mu_H \frac{d\tau}{W^2} + \frac{O(1)}{W^4}d\tau.$$

The $O(1)$ terms here may be taken to be bounded uniformly in the configuration vector $\underline{\nu}$ and matrix $\underline{\underline{3}}$, and in $W$. Choose $\varepsilon$ so that $2\mu_K - \mu_H > \varepsilon > 0$ and set $d\widetilde{\tau} = d\tau/W^2$, and use the bounds on the $O(1)$ terms to define $w^{(\varepsilon)} < \infty$ as the smallest level $w$ such that

(39) $$\left|\frac{\text{Drift } dK}{d\widetilde{\tau}} + \mu_K\right| \leq \frac{\varepsilon}{3}, \qquad \left|\frac{\text{Drift } dH}{d\widetilde{\tau}} + \mu_H\right| \leq \frac{\varepsilon}{3},$$

whenever $W \geq w$. Recall that $\ln(W) = H - 2K$, so the calculations of (38) show that $(d \ln W)^2/d\widetilde{\tau}$ is bounded, while

(40) $$\left|\frac{\text{Drift } d\ln W}{d\widetilde{\tau}} - (2\mu_K - \mu_H)\right| \leq \varepsilon$$

whenever $W \geq w^{(\varepsilon)}$. Now $\varepsilon$ was chosen so that $2\mu_K - \mu_H > \varepsilon > 0$, so it follows by consideration of the law of the iterated logarithm that if initially $W_0 > w^{(\varepsilon)}$, then there is a positive chance that $W > w^{(\varepsilon)}$ for all time; moreover, this probability increases to 1 as $W_0$ increases. In case $W > w^{(\varepsilon)}$ for all time, $W$ will grow at least linearly with rate $2\mu_K - \mu_H - \varepsilon > 0$, and hence (by considering $w^{(\varepsilon)}$ for progressively smaller $\varepsilon$)

(41) $$\frac{\ln W}{\widetilde{\tau}} \to 2\mu_K - \mu_H$$

as $\widetilde{\tau} \to \infty$.

On this event of linear growth of $W > w^{(\varepsilon)}$ the approximations in (38) improve with time. Thus as $\widetilde{\tau} \to \infty$, so the same application of the law of the iterated logarithm leads to

(42) $$\frac{K}{\widetilde{\tau}} \to -\mu_K, \qquad \frac{H}{\widetilde{\tau}} \to -\mu_H.$$

In summary, there is a positive probability of both (41) and (42) holding so long as $W_0 > w^{(\varepsilon)}$ is sufficiently large; indeed this probability increases to 1 as $W_0 \to \infty$. If $\mu_K$ and $\mu_H$ are both positive, then this ensures that $V = \exp(H)$ and $U = \exp(K)$ both hit zero (delivering coupling of both position and all stochastic areas) at $\widetilde{\tau} = \infty$.

In principle the coupling might still happen at $t$-time $\infty$, in which case it would not succeed at finite time. However,

$$d\widetilde{\tau} = \frac{d\tau}{W^2} = \left(\frac{V}{U}\right)^2 dt = \exp(2(K - H)) \, dt \tag{43}$$

and therefore the coupling will occur at $t$-time

$$\int_0^\infty \exp(-2(K - H)) \, d\widetilde{\tau}. \tag{44}$$

This will be finite on the event of linear growth of $W$ if the positive $\mu_K$ and $\mu_H$ are chosen not only to ensure that (42) holds but also so that $\mu_K < \mu_H$.

Consequently there is positive probability of coupling occurring at finite time so long as we have arranged for $\mu_K$ and $\mu_H$ to satisfy (37). □

COROLLARY 7. *The adaptive mixed coupling of Theorem 6 can be modified by adding a synchronous coupling regime so as to ensure successful coupling in finite time with probability* 1.

PROOF. We have already dealt with $n = 2$ in Section 3 above. If $n \geq 3$ and if $W$ falls below $w^{(\varepsilon)}$, so that the above procedure breaks down, then we can revert to pure synchronous coupling (33) [which holds $K$ constant and allows $H$ to evolve as a nonconstant submartingale as noted after (33)] until $W$ does exceed $w^{(\varepsilon)}$, and restart the procedure. Consequently the above can be converted into a strategy which produces coupling at finite time almost surely. □

The coupling strategy described in Corollary 7 involves discontinuous transitions between synchronous and mixture strategies, fulfilling the expectations of the heuristics at the end of Section 2. Provided we resort to time-dependent strategies, we can of course replace the mixed strategy by a time-dependent variation between reflection and rotation strategies; hence coupling can be achieved using only orthogonal controls.

**5. Complements and conclusion.** It is natural to ask whether anything might be gained by considering the full coset of coupling strategies alternate to the rotation strategies: what we might call the *rotated reflection couplings*

$$\underline{\underline{\mathbf{J}}}^{\text{rot-refl}}(\theta\underline{\underline{\mathfrak{J}}}) = (\underline{\underline{\mathbf{I}}} - 2\underline{\nu}^T\underline{\nu}) \exp(\theta\underline{\underline{\mathfrak{J}}}). \tag{45}$$

Applying the same reasoning as led to (35), we find that $\underline{\mathbf{J}}^{\text{rot-refl}}(-\tfrac{\gamma}{W}\underline{\tilde{\mathfrak{J}}})$ has the following effect:

(46)
$$(dK)^2 = \left(1 - \frac{\gamma^2}{4W^2}\|\underline{\tilde{\mathfrak{J}}}\underline{\nu}\|^2\right) d\tau + \frac{\gamma^4}{W^4}O(1)\,d\tau,$$

$$\text{Drift } dK = -\left(\frac{1}{2} - \frac{\gamma^2}{8W^2}\right) d\tau + \frac{\gamma^4}{W^4}O(1)\,d\tau,$$

$$(dK) \times (dH) = \frac{\gamma}{W}\underline{\nu}^T \underline{\underline{\mathfrak{Z}}}^T \underline{\tilde{\mathfrak{J}}}\underline{\nu}\,d\tau + \frac{\gamma^3}{W^4}O(1)\,d\tau,$$

$$(dH)^2 = 4\|\underline{\tilde{\mathfrak{J}}}\underline{\nu}\|^2 \frac{d\tau}{W^2} + \frac{\gamma^2}{W^4}O(1)\,d\tau,$$

$$\text{Drift } dH = -\left(\frac{\gamma}{2}\operatorname{tr}(\underline{\underline{\mathfrak{Z}}}^T(\underline{\mathbf{I}} - 2\underline{\nu}\,\underline{\nu}^T)\underline{\tilde{\mathfrak{J}}})\right.$$
$$\left. - (n - 1 - 4\|\underline{\underline{\mathfrak{Z}}}\underline{\nu}\|^2)\right) \frac{d\tau}{W^2} + \frac{\gamma^2 + \gamma^3}{W^4}O(1)\,d\tau.$$

This analysis would lead to a rather transparent coupling strategy if we could ensure that $H$ was always a supermartingale under a suitable rotated reflection coupling for small $\gamma/W$; however, this is not possible for $n > 3$ since it can be shown that

(47)
$$|\operatorname{tr}(\underline{\underline{\mathfrak{Z}}}^T(\underline{\mathbf{I}} - 2\underline{\nu}\,\underline{\nu}^T)\underline{\tilde{\mathfrak{J}}})| \leq \sqrt{\operatorname{tr}(\underline{\underline{\mathfrak{Z}}}_0^T \underline{\underline{\mathfrak{Z}}}_0)}$$

for $\underline{\underline{\mathfrak{Z}}}_0 = (\underline{\mathbf{I}} - \underline{\nu}\underline{\nu}^T)\underline{\underline{\mathfrak{Z}}}(\underline{\mathbf{I}} - \underline{\nu}\underline{\nu}^T)$ with the maximum being achieved when $\underline{\tilde{\mathfrak{J}}} = \underline{\underline{\mathfrak{Z}}}_0$. This maximum vanishes when $\underline{\underline{\mathfrak{Z}}}$ is of rank 2 and $\underline{\nu}$ is a nonzero eigenvector of $\underline{\underline{\mathfrak{Z}}}$, so the evolution of the configuration $(\underline{\nu}, \underline{\underline{\mathfrak{Z}}})$ unavoidably affects whether or not the drift of $H$ is negative.

It is also natural to ask whether a more direct analysis can be made using the *Carnot–Carathéodory distance* for the relevant nilpotent Lie group. Recall that the Carnot–Caratheodory distance between the origin $\underline{\mathbf{0}}$ and a point $\underline{\mathbf{x}}$ with specified stochastic areas $\underline{\underline{\mathfrak{A}}}$ is obtained by minimizing the Euclidean length of paths from $\underline{\mathbf{0}}$ to $\underline{\mathbf{x}}$ which produce the specified matrix of stochastic areas. A variational analysis shows that in general these paths are Cartesian products of circular arcs. A direct but laborious computation can be made of the stochastic calculus for the Carnot–Caratheodory distance in the two-dimensional case; unfortunately no useful picture seems to emerge from these computations.

There are various further questions to be addressed about stochastic area couplings. Certainly it is possible to use the methods described here to derive estimates on coupling rates; these are not pursued for reasons of space and also because there is a much more substantial open question:

> Can one co-adaptively couple not just the Brownian motions and their stochastic areas, but also all possible iterated path and time-integrals up to a fixed order of iteration?

Here of course it is necessary to suppose compatibility of the initial conditions, to avoid obstructions caused by algebraic relationships between the various iterated integrals (see, e.g., the algebraic remarks of Gaines [6]). Kendall and Price [13] answer this question affirmatively for the one-dimensional case by using an implicit approach; the work of this paper shows that all singly iterated path-integrals can be coupled co-adaptively, since these can all be expressed as linear combinations of Lévy stochastic areas and quadratic functions of Brownian coordinates. The general $n$-dimensional case is much more involved. We conjecture nevertheless that there is an affirmative answer to the full multidimensional question given above. However, it is clear that new approaches will have to be tried here as in the one-dimensional case: the structure which facilitates the matrix-based approach of Section 4 is no longer available for higher-order iterated integrals.

DEPARTMENT OF STATISTICS
UNIVERSITY OF WARWICK
COVENTRY CV4 7AL
UK
E-MAIL: w.s.kendall@warwick.ac.uk